\documentclass[english,12pt]{article}
\usepackage{amsfonts}
\usepackage{amssymb}
\usepackage{amssymb, amsthm, amsmath}
\usepackage[latin1]{inputenc}
\usepackage{babel}
\usepackage{graphicx}
\usepackage{psfrag}
\usepackage{amstext}
\newtheorem{teorema}{\bf Theorem}

\newtheorem{definicao}{ \bf Definition}

\newcommand{\abs}[1]{\lvert #1 \rvert}

\begin{document}
\title{ A New Characterization of the Clifford Torus}
\author{Rodrigo Ristow Montes and Jose A. Verderesi
\thanks{ristow@mat.ufpb.br and javerd@ime.usp.br}}
\date{Departamento de Matem\'atica , \\
      Universidade Federal da Para\'iba, \\[2mm]
      BR-- 58.051-900 ~ Jo\~ao Pessoa, P.B., Brazil \\
and \\
Departamento de Matem\'atica Pura, \\
Instituto de Matem\'atica e Estat\'{\i}stica, \\
Universidade de S\~ao Paulo, \\
Caixa Postal 66281, \\[2mm]
BR--05315-970~ S\~ao Paulo, S.P., Brazil}

\maketitle
\renewcommand{\thefootnote}{\fnsymbol{footnote}}
\renewcommand{\thefootnote}{\arabic{footnote}}
\setcounter{footnote}{0}
\thispagestyle{empty}
\begin{abstract}
\noindent
  In this paper we introduce the notion of contact angle for an immersed
  surface in three dimensional sphere. We deduce formulas for the Laplacian and for the
  Gaussian curvature, and we classify minimal surfaces in $S^3$ with constant
  contact angle. Also, we give an example of a minimal surface in $S^3$ with
  non  constant contact angle.
\end{abstract}
\smallskip
\noindent{\bf Keywords:} contact angle,contact distribution, Clifford
torus, minimal surfaces.
\smallskip

\noindent{\bf 2000 Math Subject Classification:} 53C42 - 53D10 - 53D35.
\section{Introduction}\mbox{}
 The notion of K\"ahler angle was introduced by  Chern and Wolfson
 in \cite{CW} and \cite{W}; it is a fundamental invariant
for  minimal surfaces in complex manifolds. Using the technique of
 moving frames, Wolfson obtained equations for the Laplacian and Gaussian
 curvature for an immersed minimal surface in $\mathbb{CP}^n$. Later, Kenmotsu
 in \cite{K}, Ohnita in \cite{O} and  Ogata in \cite{Og} classified  minimal surfaces
 with constant Gaussian curvature and constant  K\"ahler angle.\\
A few years ago, Li in  \cite{Z} gave a counterexample to the
conjecture of Bolton, Jensen and Rigoli (see \cite{BJRW}), according to which
a minimal  immersion (non-holomorphic, non anti-holomorphic, non totally real) of a  two-sphere in $\mathbb{CP}^n$
with constant  K\"ahler angle would have constant Gaussian curvature.\\
In  \cite{MV} we introduced the notion of contact angle, that can be
considered as  a new geometric invariant useful to investigate
the geometry of immersed surfaces in contact riemannian manifolds. Geometrically, the contact angle
$(\beta)$ is the complementary angle between the contact
distribution and the tangent space of the surface. Also in
\cite{MV}, we deduced formulas for the Gaussian curvature and the Laplacian
of an immersed minimal surface in $S^5$, and we found a two parameters
 family of minimal torus in $S^5$ with constant contact angle.\\
In this paper, we will construct two characterizations of the Clifford torus
 in $S^3$ using the contact angle.
We obtain the following formula for the Gaussian curvature of an immersed  minimal surface in
$S^3$:
\begin{eqnarray*}
K  &  = &  1 - \abs{\nabla\beta+e_1}^2 \nonumber
\end{eqnarray*}
Also, we obtain the following equation for the Laplacian:
\begin{eqnarray}
\Delta(\beta) &  =   & -\tan(\beta)\abs{\nabla\beta+2e_1}^2\label{eq:laplaciano1}\nonumber
\end{eqnarray}
where $e_1$ is the characteristic field defined in section 2 and  introduced by Bennequin in \cite{A}.\\
Using the equations of Gauss and  Codazzi, we have proved the following two theorems:
\begin{teorema}\label{unicidade}
The Clifford torus is the only minimal surface in $S^{3}$ with  constant contact angle.
\end{teorema}
\begin{teorema}\label{construcao}
The Clifford Torus is the only minimal surface in $S^3$ with contact angle $0
\leq \beta < \frac{\pi}{2}$ ( or $-\frac{\pi}{2} < \beta \leq 0$)
\end{teorema}

At the last section, we give two examples of minimal surfaces in $S^3$. At the
first one, we determine that the contact angle ($\beta$) of the Clifford torus
is ($\beta=0$) and the second one we determine that the contact angle of the
totally geodesic sphere is ( $\beta=arc cos(x_2)$), and therefore, non
constant (see section 6).

\section{Contact Angle for Immersed Surfaces in  $S^{3}$}\label{sec:section2}
Consider in $\mathbb{C}^{2}$ the following objects:
\begin{itemize}
\item the Hermitian product: $(z,w)= z^1\bar{w}^1+z^2\bar{w}^2$;
\item the inner product: $\langle z,w \rangle = Re (z,w)$;
\item the unit sphere: $S^{3}=\big\{z\in\mathbb{C}^{2} | (z,z)=1\big\}$;
\item the \emph{Reeb} vector field in $S^{3}$, given by: $\xi(z)=iz$;
\item the contact distribution in $S^{3}$, which is orthogonal to $\xi$:
\[\Delta_z=\big\{v\in T_zS^{3} | \langle \xi , v \rangle = 0\big\}.\]
\end{itemize}
We observe that $\Delta$ is invariant by the complex structure of $\mathbb{C}^{2}$.

Let now $S$ be an immersed orientable surface in $S^{3}$.
\begin{definicao}
The \emph{contact angle} $\beta$ is the complementary angle between the
contact distribution $\Delta$ and the tangent space $TS$ of the
surface.
\end{definicao}
Let $(e_1,e_2)$ be  a local frame of $TS$, where $e_1\in
TS\cap\Delta$. Then $\cos \beta = \langle \xi , e_2 \rangle
$.
In  $S^3$ consider the frame $f_1=z^\bot$, $f_2=iz^\bot$ and $f_3=iz$.\\
The covariant derivative  is given by:
\begin{equation}
\begin{array}{ccc}\label{eq: deriv}
Df_1  &  =  &   \;w_1^2\,f_{2}+w^{2}\,f_{3}\\
Df_2  &  =  &   \;w_2^1\,f_{1}-w^{1}\,f_{3}\\
Df_3  &  =  & -w^2\,f_{1} +w^{1}\,f_{2}
\end{array}
\end{equation}
where $(w^{1},w^{2},w^{3})$ is dual frame  of $(f_1,f_2,f_3)$.\\
Consider $e_1$ unitary vector  field in $TS\;\cap\;\Delta$, where $\Delta$ is the contact distribution.\\
Then we have:
\begin{equation}\label{eq: campos}
\begin{array}{ccl}
e_1  & = & f_{1} \\
e_2  & = & \sin(\beta)\,f_{2}+\cos(\beta)\,f_{3}\\
e_3  & = & -\cos(\beta)\,f_{2}+\sin(\beta)\,f_{3}
\end{array}
\end{equation}
where $\beta$ is the angle between $f_3$ and $e_2$,\,$(e_1,e_2)$ are tangents to $S$
and $e_3$ is normal to $S$\\
\section{Equations for the Curvature and Laplacian} \mbox{}
In this section, we deduce the equations for the Gaussian curvature and for
the Laplacian of a minimal surface in $S^3$ in terms of the contact angle.\\
Consider $(\theta^{1},\theta^{2},\theta^{3})$ dual frame  of $(e_{1},e_{2},e_{3})$
\begin{equation}
\begin{array}{ccl}\label{eq: formas}
\theta^{1}    &   =   &   w_{1}\\
\theta^{2}    &   =   &   \sin(\beta)\,w_{2}+\cos(\beta)\,w_{3}\\
\theta^{3}    &   =   &   -\cos(\beta)\,w_{2}+\sin(\beta)\,w_{3}
\end{array}
\end{equation}
At the surface $S$, we have $\theta^{3}=0$, then we obtain the equation:
\begin{equation}
\begin{array}{ccc}
sin(\beta)\,w^{3} & = & cos(\beta)w^{2}
\end{array}
\end{equation}
we have also
\begin{eqnarray}
w^2 & = & \sin \beta \theta^2 \nonumber\\
w^3 & = & \cos \beta \theta^2 \nonumber
\end{eqnarray}
It follows from (\ref{eq: formas}) that:
\begin{eqnarray}
d\theta^{1}     &   +     & \sin(\beta)(w_2^1-\cos(\beta)\theta^{2})\wedge\theta^{2}  =0\nonumber\\
d\theta^{2}     &   +     & \sin(\beta)(w_1^2+\cos(\beta)\theta^{2})\wedge\theta^{1} = 0\nonumber\\
d\theta^{3}     &   =     &   d\beta\wedge\theta^{2}-\cos(\beta)w_2^1\wedge\
w^{1}+(1+\sin^{2}(\beta))\theta^{1}\wedge\theta^{2} \nonumber
\end{eqnarray}
Therefore the connection form of $S$ is given by
\begin{eqnarray}
\theta_2^1     &   =     &   \sin(\beta)(w_2^1-\cos(\beta)\theta^{2})
\label{eq: formul}
\end{eqnarray}
Differentiating $e_3$ at the basis $(e_1,e_2)$, we have fundamental second forms coeficients
\begin{eqnarray}
De_3  &  =  &  \theta_3^1e_1+\theta_3^2e_2\nonumber
\end{eqnarray}
where
\begin{eqnarray}
\theta_3^1   &  =  & -\cos(\beta)w_2^1 - \sin^{2}(\beta)\theta^2 \nonumber\\
\theta_3^2   &  =  &  d\beta+\theta^1\nonumber
\end{eqnarray}
It follows from $d\theta^3=0$, that
\begin{eqnarray}
w_2^1(e_2)  &  =  &  -\frac {\beta_1}{\cos\beta}-\frac{(1+\sin^{2}\beta)}{\cos\beta}  \label{eq: formula1}
\end{eqnarray}
where $d\beta(e_1)=\beta_1$.\\
The condition of minimality is equivalent to the following equation
\begin{eqnarray}
\theta_1^3\wedge\theta^2 - \theta_2^3\wedge\theta^1=0\nonumber
\end{eqnarray}
we have
\begin{eqnarray}
 w_2^1(e_1) & = \frac{\beta_2}{\cos(\beta)} \label{eq: formula2}
\end{eqnarray}
where $d\beta(e_2)=\beta_2$.\\
It follows from  (\ref{eq: formul}), (\ref{eq: formula1}) and (\ref{eq:
  formula2}) that
\begin{eqnarray}
\theta_2^1  &  =  &  \tan(\beta)(\beta_2\theta^{1}-(\beta_{1}+2)\theta^{2})\nonumber \\
\theta_3^1  &  =  &  -\beta_2\theta^1+(\beta_1+1)\theta^2\nonumber\\
\theta_3^2  &  =  &   (\beta_1+1)\theta^1+\beta_2\theta^2\nonumber
\end{eqnarray}
If $J$ is the complex structure of $S$ we have $Je_1=e_2$ e $Je_2=-e_1$.\\
Using $J$, the forms above simplify to:
\begin{eqnarray}
\theta_2^1  &  =  &  \tan\beta(d\beta\circ J-2\theta^2)\nonumber \\
\theta_3^1  &  =  &  -d\beta\circ J+\theta^2\\\nonumber
\theta_3^2  &  =  &   d\beta+\theta^1\nonumber
\end{eqnarray}
It follows from Gauss equation that
\begin{eqnarray}
d\theta_1^2=\theta^1\wedge\theta^2+\theta_1^3\wedge\theta_2^3 \nonumber
\end{eqnarray}
We also have
\begin{equation}
\begin{array}{lcl}
d\theta_2^1   &   =  &  (|\nabla\beta|^2+2\beta_1)\;(\theta^{2}\wedge\theta^{1}) \label{eq: curva2}
\end{array}
\end{equation}
and therefore
\begin{eqnarray}
K  &  = &  1 - |\nabla\beta+e_1|^2 \nonumber
\end{eqnarray}Differentiating $\theta_2^1$, we have
\begin{equation}
\begin{array}{lcl}
d\theta_2^1  &  =  & \quad \sec^{2}(\beta)(|\nabla\beta|^2+2\beta_1)(\theta^{2}\wedge\theta^{1})\\
             &     & + (\tan(\beta)\Delta(\beta)+2\tan^{2}(\beta)(\beta_1+2))(\theta^{2}\wedge\theta^{1}) \label{eq: curva1}
\end{array}
\end{equation}
Using  (\ref{eq: curva2}) and (\ref{eq: curva1}), we obtain the following
formula for the Laplacian of $S$
\begin{equation}
\begin{array}{lcl}
\Delta(\beta) &  =   & -\tan(\beta)((\beta_1+2)^2+\beta^2_2)\label{eq:laplaciano1}
\end{array}
\end{equation}
Or
\begin{eqnarray}
\Delta(\beta) &  =   & -\tan(\beta)|\nabla\beta+2e_1|^2 \nonumber
\end{eqnarray}
Codazzi equations are
\begin{eqnarray}\label{eq: codazzi 1}
d\theta_1^3+\theta_2^3{\wedge}\theta_1^2=0\nonumber\\
\label{eq: codazzi 2}
d\theta_2^3+\theta_1^3{\wedge}\theta_2^1=0\nonumber
\end{eqnarray}
The first equation gives (\ref{eq:laplaciano1}) and the second equation is
always verified.
\section{Proof of the Theorem 1} \mbox{}
Suppose that  $\beta$ is constant, it follows from  (\ref{eq: curva2}) that
$d\theta_1^2=0$ and, therefore,$K=0$ ,ie., Gaussian curvature of $S$ vanishes
identically, hence $S \subset S^3$ is the Clifford torus, which prove the Theorem 1.
\section{Proof of the Theorem 2} \mbox{}
For  $0 \leq \beta < \frac{\pi}{2}$, we have $\tan \beta \geq  0$, hence
$\Delta(\beta) \leq 0$ and using that $S$ is a compact surface, we conclude
by Hopf's Lemma that $\beta$ is constant, and therefore, $K=0$ and $S$ is the
Clifford torus, which prove the Theorem 2.
\section{ Examples}
\subsection{Contact Angle of Clifford Torus in $S^3$}\mbox{}\\
Consider the torus in $S^3$ defined by:
\begin{eqnarray*}
T^2=\lbrace(z_1,z_2)\in{C^2}/z_1\bar{z_1}=\frac{1}{2},z_2\bar{z_2}=\frac{1}{2}\rbrace
\end{eqnarray*}
We consider the immersion:
\begin{eqnarray*}
f(u_1,u_2)=\frac{\sqrt 2}{2}(e^{iu_1},e^{iu_2})
\end{eqnarray*}
$T(T^2)$ is generate by $\frac{\partial}{\partial u_1}$ and
$\frac{\partial}{\partial u_2}$ it is means that:
\begin{eqnarray*}
a\frac{\partial}{\partial u_1}+b\frac{\partial}{\partial u_2}={\lambda}z^{\bot}
\end{eqnarray*}
using the condition above and the fact that $| \lambda |=1$, we obtain:
\begin{eqnarray*}
\lambda   &   =    &   ie^{i(u_1+u_2)}
\end{eqnarray*}
The unitary vector fields are:
$$\left\{
\begin{array}{lll}
e_1=ie^{i(u_1+u_2)}z^\bot\\
e_2=iz\\
e_3=e^{i\alpha}iz^{\bot}
\end{array}\right.$$
The contact angle is the angle between $e_2$ and $f_3$  ,
$$
\begin{array}{ccl}
cos(\beta) & = & \langle e_2,f_3 \rangle\\
           & = &  1
\end{array}
$$
Therefore, the contact angle is:
$$
\beta=0
$$
The fundamental second form at the basis $(e_1,e_2)$ is:\\
\begin{center}
$A=\left[\begin{array}{cc}
0 & -1\\ -1  & 0\\
\end{array}\right]$
\end{center}
\subsection{Minimal surface in $S^3$ with non constant contact angle}\mbox{}
Consider the surface described by:\\
$$\left\{
\begin{array}{lcc}
z_2-\overline{z}_2              & =  &   0\\
(x_1)^2+(y_1)^2+(x_2)^2+(y_2)^2 & =  &   1
\end{array}\right.$$\\
We see that the unitary fields are:
$$\left\{
\begin{array}{lll}
e_1=\frac{1}{\sqrt{1-x_2^{2}}}(-x_1x_2,-y_1x_2,1-x_2^{2},0)\\
e_2=\frac{1}{\sqrt{x_1^{2}+y_1^{2}}}(y_1,-x_1,0,0)\\
e_3=(0,0,0,1)
\end{array}\right.$$
The contact angle is the angle between $e_2$ and $f_3$,
$$
\begin{array}{ccl}
cos(\beta) & = & \langle e_2,f_3 \rangle\\
           & = &  x_2
\end{array}
$$
Therefore, the contact angle is:
$$
\beta= arc \, cos\, (x_2)
$$


\end{document}